\documentclass[a4paper]{paper}
\usepackage{multirow}
\usepackage{mathptmx}
\usepackage{amsmath,amsthm,amsfonts,mathtools}
\usepackage{graphicx}
\usepackage{subfig}
\DeclareMathOperator{\diag}{diag}
\DeclareMathOperator{\vec1}{vec}

\usepackage{listings}
\usepackage{subdepth} 

\usepackage{pgfplots, pgfplotstable, tikz}
\pgfplotstableset{col sep=comma}
\pgfplotsset{compat=1.16}

\usepackage{hyperref}
\newtheorem{theorem}{Theorem}

\newtheorem{remark}[theorem]{Remark}
\newtheorem{proposition}[theorem]{Proposition}
\newtheorem{definition}[theorem]{Definition}
\newcommand{\B}[1]{\mbox{\boldmath $#1$}}

\title{Relaxed Fixed Point Iterations for Matrix Equations Arising in  Markov Chains Modeling\thanks{The authors are partially supported by
    INDAM/GNCS and by the project PRA\_2020\_61 of the University of Pisa.}}

\author{
  Gemignani, Luca \\
  Dipartimento di Informatica\\
  Universit\`a di Pisa \\
  \texttt{luca.gemignani@unipi.it}
  \and
  Meini, Beatrice \\
  Dipartimento di Matematica\\
  Universit\`a di Pisa \\
  \texttt{beatrice.meini@unipi.it}
}

\begin{document}
\maketitle

\begin{abstract}
  We present some accelerated variants of fixed point  iterations for  computing the
  minimal
  non-negative solution of the unilateral matrix equation associated with  an M/G/1-type Markov chain.
  These variants derive from  certain  staircase  regular splittings of the block Hessenberg M-matrix associated with the Markov chain.   By exploiting the staircase profile we introduce
  a  two-step ﬁxed point iteration.  The iteration can   be further accelerated   by computing a weighted average between  the approximations obtained at two consecutive  steps.
  The convergence of the basic  two-step ﬁxed point iteration and of its relaxed modification is  proved. 
  Our theoretical analysis, along with  several numerical experiments, show that  the proposed variants  generally outperform the  classical iterations. 
  
\end{abstract}


\section{Introduction}
The  transition probability matrix of an M/G/1-type Markov chain  is a block Hessenberg matrix  $P$ of the form
\begin{equation}\label{maineq-1}
P=\left[\begin{array}{cccc}
    B_0 & B_1& B_2 & \ldots \\
    A_{-1} & A_0 & A_1 & \ldots \\
    &  A_{-1} & A_0 & \ddots \\
    & & \ddots & \ddots
  \end{array}\right], 
\end{equation}
with  $A_i, B_i \in \mathbb R^{n\times n}\geq 0$ and $\sum_{i=0}^\infty B_i$ and  $\sum_{i=-1}^\infty A_i$ stochastic matrices.

In the sequel, given a real matrix $A=(a_{ij})_{ij}\in\mathbb{R}^{m\times n}$, we write $A\ge 0$ ($A>0$) if $a_{ij}\ge 0$ ($a_{i,j}>0$) for any $i,j$. A stochastic matrix is matrix $A\ge 0$ such that $A\B e=\B e$, 
where $\B e$ is the column vector having all the entries equal to 1.

In the positive recurrent case, the computation of the steady state
vector $\pi$ of $P$, such that
\begin{equation}\label{maineq0}
  \B \pi^T P=\B \pi^T, \quad \B \pi^T\B e=1, \quad  \B \pi\geq \B 0,
\end{equation}
is related with the solution of the  unilateral power series matrix  equation
\begin{equation}\label{maineq}
X= A_{-1} + A_0 X + A_1 X^2 + A_2 X^3 +\ldots.
\end{equation}
Under some mild assumptions  this equation has a  componentwise minimal non-negative solution $G$ which determines, by means of Ramaswami's formula \cite{Rama},  the vector $\B \pi$.

Among the easy-to-use, but still effective, tools for  numerically solving \eqref{maineq}, there are  fixed point iterations  (see \cite{BLM}  and the references given therein for a general review of these methods).
The intrinsic simplicity of such schemes make them attractive in domains where high performance  computing is crucial.  But they come at a price:  the convergence can become very slow especially for problems
which are close to singularity.  The design of acceleration methods (aka extrapolation methods) for fixed point iterations is a classical topic in numerical analysis \cite{BRZ}. Relaxation techniques are commonly used
for  the acceleration of  classical  stationary iterative solvers for large linear systems.  In this paper we introduce some new  coupled fixed point iterations for solving  \eqref{maineq} which can be combined with
relaxation  techniques to speed up their convergence.

More specifically,  we first observe that  computing  the solution of the matrix equation   \eqref{maineq} is formally equivalent to solving a semi-infinite block Toeplitz block Hessenberg linear system.
Customary block iterative algorithms applied for the solution of this system   yield classical  fixed point iterations. In  particular,  the traditional and the U-based fixed point iteration \cite{BLM} originate
from the block Jacobi and  the block Gauss-Seidel   method, respectively.  Recently, in  \cite{GP} the authors show that some   iterative solvers  based on a block staircase partitioning  outperform
the block Jacobi and  the block Gauss-Seidel   method for M-matrix linear systems in block Hessenberg form.  The application of the staircase splitting  to the block Toeplitz block Hessenberg linear system associated with
\eqref{maineq} yields  the new coupled fixed point iteration 
\begin{equation}\label{mainit0}
  \left\{\begin{array}{ll}
      (I_n -A_0)  Y_k=A_{-1} +\sum_{i=1}^\infty A_i X_k^{i+1}; \\
      X_{k+1}= Y_k + (I_n-A_0)^{-1} A_1 (Y_k^2-X_k^2), 
  \end{array}\right. \quad k\geq 0,
  \end{equation}
starting from an initial approximation $X_0$.
 The contribution of this paper  is aimed at highlighting  the properties of \eqref{mainit0}.

We show  that, if $X_0=0$, the sequence $\{ X_k\}_k$ defined in  \eqref{mainit0} converges  to $G$ faster than the traditional fixed point iteration.
In the case where the starting matrix $X_0$  of  \eqref{mainit0}  is any column stochastic matrix and $G$ is also column stochastic,  we prove that
  the sequence $\{X_k\}_{k}$
  still  converges  to $G$.  Moreover, by comparison of the  mean asymptotic rates  of
  convergence, we conclude that \eqref{mainit0} is  asymptotically faster than  the  traditional fixed point iteration.

At each iteration the scheme \eqref{mainit0} determines two approximations which   can be combined by
using the relaxation technique, that is, the approximation computed at the $k$-th step takes the form  of   a weighted average   between  $Y_k$ and  $X_{k+1}$.  The   modified  relaxed variant
is defined by the sequence  \begin{equation}\label{mainitrel0}
  \left\{\begin{array}{ll}
      (I_n -A_0)  Y_k=A_{-1} +\sum_{i=1}^\infty A_i X_k^{i+1}; \\
      X_{k+1}= Y_k + \omega_{k+1} (I_n-A_0)^{-1} A_1 (Y_k^2-X_k^2), 
  \end{array}\right. \quad k\geq 0,
  \end{equation}
where $\omega_{k+1}$ is the relaxation parameter.
The convergence results for  \eqref{mainit0} easily extend to  the modified  scheme in  the case of under-relaxation, that is,  the parameter $\omega_k$ satisfies $0\leq \omega_k \leq 1$.
Heuristically,  it is argued  that  over-relaxation values ($\omega_k >1)$  can improve  the  convergence.   If $X_0=0$, under some assumptions a theoretical estimate of
the asymptotic convergence rate of \eqref{mainitrel0} is given which confirms this  heuristic. Moreover,  an adaptive  strategy is devised which makes possible  to perform over-relaxed
iterations of  \eqref{mainitrel0} by still ensuring the convergence of the  overall iterative process.  
   The results of extensive numerical  experiments  confirm the  effectiveness of the proposed  variants,
  which  generally outperform  the $U$-based fixed point iteration for nearly singular problems. In particular,     the over-relaxed scheme
  \eqref{mainitrel0}  with  $X_0=0$, combined  with the adaptive strategy for parameter estimation,  is capable to significantly accelerate  the convergence without increasing the computational cost.

  The paper is organized as follows. In Section~\ref{SEC2} we set up the theoretical framework, briefly recalling some preliminary properties and assumptions.
  In Section~\ref{SEC3}  we  revisit classical fixed point iterations for solving \eqref{maineq}  by establishing the link with the iterative solution of an associated
  block Toeplitz block Hessenberg linear system. In Section~\ref{SEC4}   we  introduce the  new fixed point iteration \eqref{mainit0}  by proving some  convergence results.
  The  relaxed variant  \eqref{mainitrel0},  as well as the  generalizations of convergence results for  this variant, are described in Section~\ref{SEC5}. Section~\ref{SEC6} deals with
  a formal analysis of the asymptotic convergence rate of both  \eqref{mainit0} and \eqref{mainitrel0}.  Adaptive strategies for the  choice of the
  relaxation parameter  are discussed in Section~\ref{SEC7}, together with their cost analysis under some simplified assumptions.   Finally, the results of extensive numerical  experiments   are presented in
  Section~\ref{SEC8} whereas  conclusions and future work are  the subjects   of Section~\ref{SEC9}.

\section{Preliminaries and assumptions}\label{SEC2}
Throughout the paper we assume that $A_i$, $i\geq -1$, are $n\times n$ nonnegative matrices, such that their sum
 $A=\sum_{i=-1}^{\infty} A_i$ is irreducible and  row stochastic, that is, $A\B e=\B e$, $\B e=\left[1, \ldots, 1\right]^T$.
  According to the results of \cite[Chapter 4]{BLM}, such assumption implies that  \eqref{maineq} has a unique  componentwise minimal nonnegative solution $G$; moreover, 
  $I_n -A_0$ is an invertible  M-matrix and, hence, $(I-A_0)^{-1}\geq 0$.   
  
  Furthermore, in view of the Perron Frobenius Theorem,
  there exists a  unique vector $\B v$ such that $\B v^T A=\B v^T$ and $\B v^T \B e=1$, $\B v>0$.   
  If the series $\sum_{i=-1}^{\infty} iA_i$ is convergent we may define the vector $\B w=\sum_{i=-1}^{\infty} iA_i\B e\in \mathbb { R}^n$.
 In the study of M/G/1-type Markov chains, the {\em drift} is the scalar number  $\eta=\B v^T \B w$ \cite{NMF}. The sign of the drift determines the positive recurrence of the M/G/1-type Markov chain \cite{BLM}.
  
  When explicitly stated, we will assume the following condition:
  \begin{enumerate}
  \item[A1.]\label{it3}   The series $\sum_{i=-1}^{\infty} iA_i$ is convergent and $\eta<0$.  
  \end{enumerate}
 Under assumption [A1],  
 the componentwise minimal nonnegative  solution $G$ of equation \eqref{maineq} is stochastic, i.e.,
   $G\B e= \B e$ (\cite{BLM}). Moreover $G$ is the only stochastic solution.

\section{Nonlinear matrix equations and structured linear systems}\label{SEC3}

In this section we reinterpret classical fixed point iterations for solving the matrix equation \eqref{maineq} in terms of iterative methods for solving a structured linear system.

Formally,  the power series matrix equation \eqref{maineq} can be rewritten as the following block Toeplitz block Hessenberg   linear system
\[
\left[\begin{array}{cccc}
    I_n-A_0 & -A_1& -A_2 & \ldots \\
    -A_{-1} & I_n -A_0& -A_1 & \ldots \\
    & - A_{-1} & I_n-A_0 & \ddots \\
    & & \ddots & \ddots
  \end{array}\right]\left[\begin{array}{cccc}
    X \\
    X^2  \\
    X^3 \\
    \vdots
  \end{array}\right]=\left[\begin{array}{cccc}
    A_{-1} \\
    0  \\
    0 \\
    \vdots
  \end{array}\right].
\]
The above linear system  can  be expressed in compact form as
\begin{equation}\label{fixeq0}
 H \B {\hat X} = \B E A_{-1},  \  H=I -\tilde P, 
\end{equation}
where $\B {\hat X}^T=\left[X^T, {X^T}^2, \ldots\right]^T$, $\tilde P$  is the matrix obtained from $P$ in \eqref{maineq-1} by  removing its  first block row and block column and
$\B E=\left[I_n, 0_n, \ldots\right]^T$.
Classical fixed point iterations for solving  \eqref{maineq} can be interpreted as iterative methods for solving \eqref{fixeq0}, based on suitable partitionings of the matrix $H$.

For instance, from the partitioning $H=M-N$, where $M=I$ and $N=\tilde P$,
we find that
the block vector $\hat X$ is a solution of the fixed point problem
\begin{equation}\label{fixeq}
\B {\hat X}=\tilde P \B {\hat X} +  \B E A_{-1}.
\end{equation}
From this equation we may generate the sequence of block vectors
\[
\B {\hat X_k}=
\begin{bmatrix}
X_k\\
X_k^2\\
X_k^2\\
\vdots
\end{bmatrix},~~~
\B {Z_{k+1}}=
\begin{bmatrix}
X_{k+1}\\
X_{k+1} X_{k}\\
X_{k+1} X_{k}^2\\
\vdots
\end{bmatrix},
\]
such that
\[
\B Z_{k+1}=\tilde P \B {\hat X_k} +  \B E A_{-1},~~k=0,1,\ldots.
\]
We may easily verify that the sequence $\{X_k\}_k$ coincides with the sequence generated by the so called  {\em natural fixed point iteration} $X_{k+1}=\sum_{i=-1}^\infty A_i X_k^{i+1}$, $k=0,1,,\ldots$, applied to \eqref{maineq}.

Similarly, the 
 {\em Jacobi partitioning}, where $M=I\otimes (I_n-A_0)$ and $N=M-H$,
 which leads to the sequence
 \[
M \B Z_{k+1}= N \B {\hat X_k} +  \B E A_{-1},~~k=0,1,\ldots,
\]
corresponds to the {\em traditional fixed point iteration}
\begin{equation}\label{eq:trad}
(I_n-A_0)X_{k+1}= A_{-1} +\sum_{i=1}^\infty A_i X_k^{i+1},~~k\ge0.
\end{equation}
The {\em anti-Gauss-Seidel  partitioning},   where  $M$ is the block  upper triangular part of
$H$ and $N=M-H$, determines the {\em $U$-based  fixed point iteration}
\begin{equation}\label{eq:ubased}
\left(I_n-\sum_{i=0}^\infty A_iX_k^i\right)X_{k+1}=A_{-1},~~k\ge 0.
\end{equation}
The convergence properties of these three fixed point iterations are analyzed in  \cite{BLM}. Among  the three iterations  \eqref{eq:ubased} is the fastest  and also the most expensive since
it requires the solution of a new linear system (with multiple right hand sides) at each iteration.  Moreover,  
it turns out that fixed-point iterations  exhibit  arbitrarily slow convergence for problems which are close to singularity.  In particular, for
positive recurrent Markov chains having a drift  close to zero the convergence slows
down and the number of iterations  becomes arbitrarily large.  In the next sections  we present some new  fixed point iterations   which offer several advantages in terms
of computational efﬁciency and convergence properties when compared with \eqref{eq:ubased}.

\section{A new fixed point iteration}\label{SEC4}
Recently in \cite{GP}  a comparative analysis has been performed for
the asymptotic convergence rates of some regular splittings  of a non-singular  block upper Hessenberg  M-matrix of finite size.  The conclusion is that  the  staircase splitting is faster than the
anti-Gauss-Seidel splitting, that  in turn is faster  than the Jacobi splitting.  The second result is classical, while the first one is  somehow surprising since  the  matrix $M$  in the staircase
splitting is  much more sparse than the corresponding matrix  in the anti-Gauss-Seidel partitioning and the  splittings are not comparable. 
Inspired from these convergence properties, we introduce a new fixed point iteration for solving \eqref{maineq}, based on the {\em staircase  partitioning} of $H$, namely, 
\begin{equation}\label{mstair}
M=\left[\begin{array}{cccccccccc}
    I_n-A_0  \\
    -A_{-1} & I_n-A_0 & -A_1 \\
    &   & I_n-A_0  \\
    & & -A_{-1} & I_n-A_0 & -A_1\\
    &&&&  I_n-A_0 \\
    &&&&\times & \times  & \times  \\
    &&&&&&\times &\phantom{a}
    \\
    \\
  \end{array}\right], \quad N=M-H.
\end{equation}

The splitting has attracted  interest for applications in  parallel computing environments \cite{Meu,S1}.  
In principle the alternating  structure of the matrix $M$ in \eqref{mstair} suggests several  different iterative schemes.

From one hand,  the odd block entries of the system $M \B Z_{k+1}=N\B {\hat X_k} +  \B E A_{-1}$
  yield the traditional fixed point iteration.
On the other hand, the even block entries
lead  to the  implicit scheme $-A_{-1} + (I_n-A_0)X_{k+1}  -A_1 X_{k+1}^2 =\sum_{i=2}^\infty A_iX_k^{i+1}$   recently  introduced in \cite{BLMfamily}.
Differently,  by looking at the structure of the matrix $M$ on the whole,  we introduce the following  composite two-stage iteration:
\begin{equation}\label{eq:0}
  \left\{\begin{array}{ll}
      (I_n -A_0)  Y_k=A_{-1} +\sum_{i=1}^\infty A_i X_k^{i+1}; \\
      -A_{-1} + (I_n-A_0)X_{k+1}  -A_1 Y_k^2=\sum_{i=2}^\infty A_iX_k^{i+1},
  \end{array}\right. \quad k\geq 0,
  \end{equation}
  or, equivalently,
  \begin{equation}\label{mainit}
  \left\{\begin{array}{ll}
      (I_n -A_0)  Y_k=A_{-1} +\sum_{i=1}^\infty A_i X_k^{i+1}; \\
      X_{k+1}= Y_k + (I_n-A_0)^{-1} A_1 (Y_k^2-X_k^2), 
  \end{array}\right. \quad k\geq 0,
  \end{equation}
  starting from an initial approximation $X_0$.
  At each step $k$, this scheme consists of a traditional fixed point iteration that computes $Y_k$ from $X_k$, followed by a cheap correction step for computing the new approximation $X_{k+1}$.

As  for classical fixed point iterations,  the  convergence is guaranteed  when   $X_0=0$.
  \begin{proposition}\label{prop1}
    Assume that $X_0=0$. Then the sequence $\{X_k\}_{k\in \mathbb N}$ generated by \eqref{mainit} converges monotonically to $G$.
  \end{proposition}
  \begin{proof}
    We show by induction on $k$ that $0\le X_k\le Y_k\le X_{k+1}\le G$ for any $k\ge 0$.  
    For $k=0$ we verify easily that 
    \[
    X_1\geq Y_0=(I_n-A_0)^{-1}A_{-1}\geq 0=X_0, \quad  (I_n -A_0)X_1\leq A_{-1} +A_1G^2\leq (I_n -A_0)G, 
    \]
    which gives $G\geq X_1$.
    Suppose  now that $G\geq X_k\geq Y_{k-1}\geq X_{k-1}$, $k\ge 1$.  We find that
  \[
  (I_n -A_0)  Y_k=A_{-1} +\sum_{i=1}^\infty A_i X_k^{i+1}\geq A_{-1} + A_1 Y_{k-1}^2 + \sum_{i=2}^\infty A_iX_{k-1}^{i+1} =(I_n -A_0)  X_k.
  \]
  and
  \[
  (I_n -A_0)  Y_k=A_{-1} +\sum_{i=1}^\infty A_i X_k^{i+1}\leq A_{-1} +\sum_{i=1}^\infty A_i G^{i+1}=(I_n -A_0)  G.
  \]
  By multiplying both sides by the inverse of $I -A_0$  we obtain  that  $G\geq Y_k\geq X_k$.  This also implies that $Y_k^2-X_k^2\geq 0$ and therefore $X_{k+1}\geq  Y_k$.
  Since
  \[
  (I_n -A_0)  X_{k+1}=A_{-1} + A_1 Y_k^2+\sum_{i=2}^\infty A_i X_k^{i+1}\leq (I_n -A_0)  G,
  \]
  we  prove  similarly that $G\geq X_{k+1}$. It follows that $\{X_k\}_{k\in \mathbb N}$ is convergent, the limit  solves  \eqref{maineq} by continuity and, hence, the limit coincides with the matrix $G$, since $G$ is the minimal nonnegative solution.
  \end{proof}

 A similar  result also holds  for the case where $X_0$ is a stochastic matrix, assuming that [A1] holds, so that $G$ is  also stochastic.
  
 . 
  \begin{proposition}\label{propsto}
    Assume that condition [A1] is  fulfilled and that $X_0$  is a stochastic matrix.
    Then, the sequence $\{X_k\}_{k\in \mathbb N}$ generated by \eqref{mainit} converges to $G$.
  \end{proposition}
  \begin{proof}
  From \eqref{eq:0}  we obtain that 
  \[
  \left\{\begin{array}{ll}
      (I_n -A_0)  Y_k=A_{-1} +\sum_{i=1}^\infty A_i X_k^{i+1}; \\
     (I_n -A_0)  X_{k+1}=A_{-1} + A_1 Y_k^2+\sum_{i=2}^\infty A_i X_k^{i+1}
  \end{array}\right. \quad k\geq 0,
   \]
 which  gives that  $X_k\geq 0$ and $Y_k\ge 0$, for any $k\in \mathbb N$, since $X_0\ge 0$. By assuming that $X_0\B e=\B e$, we may easily show by induction that $Y_k\B e= X_k\B e=\B e$ for any $k\ge 0$. Therefore, all the matrices  $X_k$ and $Y_k$, $k\in \mathbb N$, are stochastic.
Let $\{\hat X_k\}_{k\in \mathbb N}$ be the sequence generated by \eqref{mainit}  with   $\hat X_0=0$. We can easily  show by induction that $X_k\geq \hat X_k$ for any $k\in \mathbb N$. 
Since $\lim_{k\to\infty}\hat X_k=G$, then any convergent subsequence of $\{X_k\}_{k\in \mathbb N}$  converges to  a stochastic matrix $S$ such that $S\geq G$.
Since $G$ is also stochastic,   it follows that $S=G$ and therefore, by compactness,
we conclude that   the sequence $\{X_k\}_{k\in \mathbb N}$ is   also convergent to $G$. 
  \end{proof}

  Propositions \ref{prop1} and \ref{propsto} are global convergence  results. An estimate of the rate of convergence of \eqref{mainit} will be
  provided  in Section~\ref{SEC6}, 
  together  with a comparison with  other existing methods.

  \section{A relaxed variant}\label{SEC5}

  At each iteration, the scheme \eqref{mainit} determines two approximations which   can be combined by
  using a relaxation technique, that is, the approximation computed at the $k$-th step takes the form  of   a weighted average   between  $Y_k$ and  $X_{k+1}$,
   \[
  X_{k+1}= \omega_{k+1} (Y_k + (I_n-A_0)^{-1} A_1 (Y_k^2-X_k^2)) + (1-\omega_{k+1})Y_k, \quad k\geq 0, 
  \]
 In matrix  terms, the resulting relaxed variant of   \eqref{mainit} can be written as 
  \begin{equation}\label{mainitrel}
  \left\{\begin{array}{ll}
      (I_n -A_0)  Y_k=A_{-1} +\sum_{i=1}^\infty A_i X_k^{i+1}; \\
      X_{k+1}= Y_k + \omega_{k+1} (I_n-A_0)^{-1} A_1 (Y_k^2-X_k^2), 
  \end{array}\right. \quad k\geq 0.
  \end{equation}
   If $\omega_k=0$, for $k\geq 1$,  the relaxed scheme reduces to the  traditional fixed point iteration \eqref{eq:trad}. If $\omega_k=1$, for $k\geq 1$,  the relaxed scheme  coincides with  \eqref{mainit}.
  Values of $\omega_k$ greater than 1 can speed up the convergence of the iterative scheme.  
  
  Concerning convergence,
  the proof of Proposition \ref{prop1}  can immediately  be generalized to show that the sequence $\{ X_k\}_k$ defined  by \eqref{mainitrel},  with   $X_0=0$, converges for any  $\omega_{k}=\omega$, $k\geq 1$,  such that $0\leq \omega\leq 1$.  Moreover,
      let  $\{X_k\}_{k\in \mathbb N}$ and $\{\hat X_k\}_{k\in\mathbb N}$ be the sequences generated by \eqref{mainitrel} for $\omega_{k}=\omega$ and $\omega_{k}=\hat \omega$ with $0\leq \omega\leq \hat \omega\leq 1$,
          respectively.  It can be easily shown that $G\geq \hat X_k \geq X_k$  for any $k$ and, hence, 
          that the iterative scheme \eqref{mainit} converges faster than  \eqref{mainitrel} if $0\leq \omega_{k}=\omega<1$.

  The convergence analysis of the modified scheme \eqref{mainitrel} for $\omega_k>1$  is  much more involved since   the choice of a  relaxation parameter $\omega_k>1 $ can destroy the monotonicity  and the
  nonnegativity of the approximation sequence, 
which is at the core of  the proofs  of  Proposition \ref{prop1} and Proposition \ref{propsto} .
In order to maintain the convergence properties  of the modified scheme  we  introduce  the following definition.

\begin{definition}\label{deff}
  The sequence $\{\omega_k\}_{k\geq 1}$ is {\em eligible}  for the scheme \eqref{mainitrel}  if $\omega_k\geq 0$, $k\geq 1$,  and the following two conditions are satisfied:
  \begin{equation}\label{eqom1}
  \omega_{k+1}A_1(Y_k^2-X_k^2)\leq A_1(X_{k+1}^2-X_k^2) +\sum_{i=2}^{\infty} A_i(Y_k^{i+1}-X_k^{i+1}), \quad k\geq 0.
  \end{equation}
  and
  \begin{equation}\label{eqom2}
   X_{k+1}\B e =Y_k\B e + \omega_{k+1} (I_n-A_0)^{-1} A_1 (Y_k^2-X_k^2)\B e\leq \B e, \quad k\geq 0.
  \end{equation}
\end{definition}
It is worth noting  that  condition \eqref{eqom1} is  implicit since the construction of $X_{k+1}$  also depends on   the value of $\omega_{k+1}$.  By replacing  $X_{k+1}$ in \eqref{eqom1} with the  expression in
the right-hand side  of  \eqref{mainitrel} we obtain  a quadratic inequality  with matrix coefficients in the variable  $\omega_{k+1}$.
Obviously $\omega_k=\omega$, $0\leq \omega \leq 1$, $k\geq 1$, are  eligible  sequences.

The following  generalization of Proposition \ref{prop1} holds.
\begin{proposition}\label{prop2}
  Set $X_0=0$ and let condition [A1] be satisfied. If   $\{\omega_k\}_{k\geq 1}$ is eligible  then  the sequence $\{X_k\}_{k\in \mathbb N} $ generated by \eqref{mainitrel} converges monotonically to $G$. 
\end{proposition}
\begin{proof}
 We show by induction that  $0\le X_k \le Y_k\le X_{k+1}\le G$.  For $k=0$ we have  
 \[
 (I_n-A_0)X_1\geq (I_n-A_0)Y_{0}=A_{-1}\geq 0=X_0
 \]
 which gives immediately  $X_1\geq Y_{0}\geq 0$. Moreover, $X_1\B e\le \B e$.
    Suppose now  that $ X_k\geq Y_{k-1}\geq X_{k-1}\geq 0$, $k\ge 1$.  We find that
    \[
    \begin{array}{lll}
      (I_n  -A_0)  X_{k}=A_{-1} + A_1(X_{k-1}^2 +\omega_k (Y_{k-1}^2 -X_{k-1}^2))+ \sum_{i=2}^\infty A_i X_{k-1}^{i+1}\leq \\
      \leq A_{-1}+A_1 X_{k-1}^2 + A_1 (X_k^2 -X_{k-1}^2) +  \sum_{i=2}^\infty A_i(Y_{k-1}^{i+1}-X_{k-1}^2)+ \sum_{i=2}^\infty A_i X_{k-1}^{i+1}\leq\\
      \leq  A_{-1} + A_1 X_k^2 + \sum_{i=2}^\infty A_i Y_{k-1}^{i+1}\leq  A_{-1} + \sum_{i=1}^\infty A_i X_{k}^{i+1}=(I_n -A_0)  Y_k
      \end{array}
  \]  
  from which it follows $ Y_k\geq X_k\geq 0$.  This also implies that $X_{k+1}\geq Y_k$. From \eqref{eqom2} it follows that $X_k \B e\leq \B e$ for all $k\geq 0$ and therefore
  the sequence of approximations is upper bounded and it has a finite limit $H$.  By continuity   we  find that $H$ solves the matrix equation \eqref{maineq} and $H\B e\le \B e$. Since $G$ is the unique stochastic solution, then $H=G$.
\end{proof}


\begin{remark}\label{rk2}
   As previously mentioned,  condition \eqref{eqom1} is implicit, since  $X_{k+1}$ also depends on $\omega_{k+1}$.  An explicit condition can be derived by noting that
  \[
  A_1(X_{k+1}^2-X_k^2) \geq A_1((Y_k^2-X_k^2) +\omega_{k+1}( Y_k\Gamma_k + \Gamma_k Y_k)),
  \]
  with $\Gamma_k=(I_n-A_0)^{-1} A_1 (Y_k^2-X_k^2)$.
  There follows that \eqref{eqom1} is fulfilled whenever
  \[
    \frac{\omega_{k+1}-1}{\omega_{k+1}} A_1((Y_k^2-X_k^2)\leq ( Y_k\Gamma_k + \Gamma_k Y_k) + \omega_{k+1}^{-1} \sum_{i=2}^{\infty} A_i(Y_k^{i+1}-X_k^{i+1})
    \]
    which  can be reduced to a linear inequality in $\omega_{k+1}$  over a fixed search interval. 
    Let $\omega\in [1, \hat \omega]$  be such that 
    \begin{equation}\label{strat1}
    \frac{\omega_{}-1}{\omega_{}} A_1((Y_k^2-X_k^2)\leq ( Y_k\Gamma_k + \Gamma_k Y_k) + {\hat \omega}^{-1} \sum_{i=2}^{\infty} A_i(Y_k^{i+1}-X_k^{i+1}).
    \end{equation}
    Then we can impose that 
    \begin{equation}\label{strat2}
      \omega_{k+1}=\max\{\omega_{}\colon \omega \in  [1, \hat \omega]  \land  \eqref{strat1} \ holds\}.
    \end{equation}
    From a computational viewpoint   the strategy based on   \eqref{strat1}  and \eqref{strat2} for the choice of the value of  $\omega_{k+1}$
    can be  too much expensive and some weakened criterion   should be considered (compare with Section~\ref{SEC7} below). 
\end{remark}   

In the following section we perform a convergence analysis to estimate the convergence rate of \eqref{mainitrel} in the stationary case $\omega_k=\omega$, $k\geq 1$, as a function of the relaxation parameter.

\section{Estimate of the convergence rate}\label{SEC6}
   
    Relaxation techniques are
   usually aimed at accelerating the convergence  speed of  frustratingly  slow  iterative solvers. Such  inefficient behavior is  typically  exhibited
   when the solver is applied to a nearly singular problem. Incorporating some relaxation parameter into the iterative scheme \eqref{maineq} can greatly improve its convergence rate. 
  Preliminary insights  on the effectiveness of relaxation techniques applied for the solution of the fixed point problem  \eqref{fixeq} come from  the classical analysis for stationary iterative solvers and are developed in Section \ref{sec:fin}. A precise convergence analysis is presented in Section \ref{sec:conv}.

  \subsection{Finite dimensional convergence analysis}\label{sec:fin}
\begin{figure}
  \includegraphics[width=0.5\textwidth]{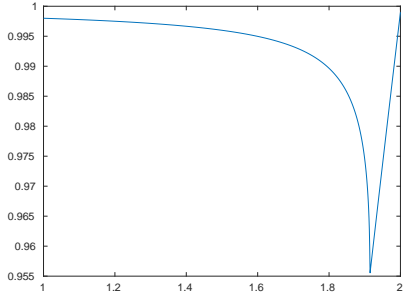}
  \caption{Plot of $\rho_S(\omega)$ for $\omega\in [1,2]$ and $\lambda=0.999$}
\label{fig1}
\end{figure}

  Suppose that $H$   in \eqref{fixeq} is block tridiagonal of finite size $m=n\ell$, $\ell$ even.  We are interested in comparing the iterative algorithm based on the splitting \eqref{mstair} with other classical iterative solvers for the
  solution of a linear system with coefficient matrix $H$.  As usual, we can write
  $H=D -P_1 -P_2$, where $D$ is block diagonal, while $P_1$ and $P_2$ are staircase matrices with zero block diagonal.  The eigenvalues $\lambda_i$ of  the Jacobi   iteration matrix satisfy
  \[
 0=\det(\lambda_i I_\ell -D^{-1}P_1 -D^{-1}P_2).
 \]
 Let us consider a relaxed  scheme  where  the matrix $M$   is obtained from \eqref{mstair} by multiplying the off-diagonal blocks  by $\omega$.
  The eigenvalues $\mu_i$ of the iteration matrix associated with the relaxed  staircase regular splitting  are such that
  \[
  0=\det(\mu_i (D-\omega P_1) -(P_2+(1-\omega) P_1))
  \]
  and, equivalently,
  \[
  0=\det(\mu_i I_\ell -(\mu_i \omega +(1-\omega))D^{-1}P_1 -D^{-1}P_2).
  \]
  By  using a  similarity transformation induced by the matrix $S=I_{\ell/2} \otimes  \diag\left[I_n,\alpha I_n\right]$  we find that
  \begin{align*}
  &\det(\mu_i I -(\mu_i \omega +(1-\omega))D^{-1}P_1 -D^{-1}P_2)= \\
  &\det(\mu_i I_\ell -\alpha (\mu_i \omega +(1-\omega))D^{-1}P_1 -\frac{1}{\alpha}D^{-1}P_2).
  \end{align*}
  There follows that
  \[
  \mu_i\alpha=\lambda_i
  \]
  whenever $\alpha$ fulfills
  \[
  \alpha (\mu_i \omega +(1-\omega))=\frac{1}{\alpha}.
  \]
  Therefore,  the eigenvalues of the Jacobi  and  relaxed  staircase regular splittings are related by
  \[
  \mu_i^2 -\lambda_i^2 \mu \omega + \lambda_i^2 (\omega-1)=0.
  \]
  For $\omega=0$ the staircase splitting reduces to the Jacobi partitioning.
  For $\omega=1$ we find that $\mu_i=\lambda_i^2$  which yields the classical relation between the spectral radii of  Jacobi and  Gauss-Seidel methods. Observe that  it is well
  known that the asymptotic convergence rates of Gauss-Seidel and  the  staircase iteration coincide  when applied to a block tridiagonal matrix \cite{AM}.
  For $\omega>1$  the spectral radius of the relaxed staircase  scheme can be significantly less than the  spectral radius of  the  same scheme for $\omega=1$.
  In Figure \ref{fig1} we  illustrate the plot of the function
  \[
  \rho_S(\omega)=\left|\frac{\lambda^2 \omega + \sqrt{\lambda^4 \omega^2 -4 \lambda^2 (\omega-1)}}{2}\right|, \quad 1\leq \omega \leq 2, 
    \]
  for a fixed $\lambda=0.999$. 
For the best choice  of $\omega=\omega^\star=2 \displaystyle\frac{ 1 + \sqrt{1-\lambda^2}}{\lambda^2}$   we find
$\rho_S(\omega^\star)=1-\sqrt{1-\lambda^2} =\displaystyle\frac{\lambda^2}{1+\sqrt{1-\lambda^2}}$.

\subsection{Asymptotic Convergence Rate}\label{sec:conv}
A formal analysis of the asymptotic convergence rate of  the relaxed  variants  \eqref{mainitrel} can be carried out  by using the tools  described in \cite{MeiB}.  In this section we relate the approximation error at two subsequent steps and we provide an estimate of the asymptotic rate of convergence, expressed as the spectral radius of a suitable matrix depending on $\omega$.

Hereafter it is assumed that assumption [A1] is verified. 

\subsubsection{The case  $X_0=0$}
Let us introduce the error matrix $E_k=G -X_k$, where $\{X_k\}_{k\in \mathbb N}$ is generated by \eqref{mainitrel} with
$X_0=0$. We also define $E_{k+1/2}=G-Y_k$, for $k=0,1,2,\ldots$.
Suppose that
\begin{enumerate}
  \item[C0.]\label{itc0}    $\{\omega_k\}_k$
   is an eligible  sequence according to Definition \ref{deff}.  
\end{enumerate}
Under this assumption from  Proposition \ref{prop2} the sequence $\{X_k\}_k$ converges monotonically to $G$ and $E_k\ge 0$, $E_{k+1/2}\ge 0$.
Since $E_k\geq 0$ and $\parallel E_k\parallel_\infty=\parallel E_k\B e\parallel_\infty$,  we analyze the convergence of
the vector $\B \epsilon_k= E_k\B e$, $k\geq 0$.

We have 
\begin{equation}\label{err1}
(I_n-A_0) E_{k+1/2}= \sum_{i=1}^\infty A_i (G^{i+1}-X_k^{i+1})= \sum_{i=1}^\infty A_i\sum_{j=0}^i G^j E_kX_{k}^{i-j}.
\end{equation}
Similarly, for the second equation of \eqref{mainitrel}, we find that
\[
(I_n-A_0) E_{k+1}= (I_n-A_0) E_{k+1/2} -\omega_{k+1} A_1((G^2-X_k^2)-(G^2-Y_k^2)),
\]
which gives
\begin{align}\label{err2}
  &(I_n-A_0) E_{k+1}=  \nonumber\\
  &(I_n-A_0) E_{k+1/2}- \omega_{k+1}A_1(E_kG+X_k E_k) +\omega_{k+1}A_1(GE_{k+1/2}+E_{k+1/2} Y_k).
\end{align}
Denote by $R_k$ the matrix on the right hand side of \eqref{err1}, i.e., 
\[R_k=\sum_{i=1}^\infty A_i\sum_{j=0}^i G^j E_kX_{k}^{i-j}.
\]
Since $G\B e=\B e$, equation \eqref{err2}, together with the monotonicity, yields 
\begin{equation}\label{err3}
  (I_n-A_0) E_{k+1}\B e\le  R_{k}\B e- \omega_{k+1}A_1(I_n+X_k) E_k\B e +\omega_{k+1}A_1(I_n+G)(I_n-A_0)^{-1}R_k \B e.
  \end{equation}
  Observe that $R_k \B e\le  W E_k \B e$, where
  \[
  W=\sum_{i=1}^\infty A_i\sum_{j=0}^i G^j,
  \]
 hence
  \begin{equation}\label{asyr}
  \epsilon_{k+1}\le P(\omega_{k+1})\epsilon_k, \quad k\geq 0,
  \end{equation}
  where
  \begin{equation}\label{defpomega}
  P(\omega)=
  (I_n-A_0)^{-1}W - \omega(I_n-A_0)^{-1} A_1(I_n+G)(I_n-(I_n-A_0)^{-1}W).
  \end{equation}
  The matrix $P(\omega)$ can be written as $P(\omega)=M^{-1}N(\omega)$ where 
  \[
  M=I_n-A_0,~~~ M-N(\omega)=A(\omega),
  \]
  and
  \[
  A(\omega)=(I_n+\omega A_1(I_n+G)(I_n-A_0)^{-1})(I_n-A_0-W).
  \]
  Let us assume the  following condition holds: 
  \begin{enumerate}
  \item[C1.]\label{itc1}   The relaxation  parameter $\omega$  satisfies  $\omega\in [0,\hat\omega]$ with
   $\hat \omega\ge 1 $ such that 
  \[
\hat \omega  A_1(I_n+G)(I_n-A_0)^{-1}(I_n-A_0-W)\leq W.
\]
\end{enumerate}
  Assumption [C1]  ensures that $N(\omega)\geq 0$ and, therefore $P(\omega)\geq 0$ and $M-N(\omega)=A(\omega)$ is a regular splitting of $A(\omega)$.  If
  [C1]  is satisfied at  each iteration of \eqref{mainitrel}, then 
  from  \eqref{asyr} we obtain that
  \[
  \epsilon_k\le P(\omega_k)P(\omega_{k-1})\cdots P(\omega_0)\epsilon_0,~~k\ge 1.
  \]
  Therefore, the asymptotic rate of convergence, defined as
  \[
  \sigma=\limsup_{k\to\infty}\left( \frac{\| \B \epsilon_k\|}{ \|\B  \epsilon_0 \|}\right)^{1/k},
  \]
  where $\|\cdot\|$ is any vector norm, is such that
  \[
  \sigma\le \limsup_{k\to\infty}\| P(\omega_k)P(\omega_{k-1})\cdots P(\omega_0)\|_\infty^{1/k}.
  \]
The above properties can be summarized in the following result, that
  gives a  convergence rate estimate for iteration \eqref{mainitrel}. 

  \begin{proposition}\label{PF0}
    Under Assumptions [A1], [C0] and [C1],  for the fixed point iteration \eqref{mainitrel}  applied  with $\omega_k=\omega$ for any $k\ge 0$, we have the following
    convergence rate estimate:
    \[
    \sigma=\limsup_{k\to\infty}\left( \frac{\|\B \epsilon_k\|}{ \| \B \epsilon_0 \|}\right)^{1/k}\leq \rho_\omega,
    \]
    where $P(\omega)$ is defined in \eqref{defpomega} and $\rho_\omega=\rho(P(\omega))$ is the spectral radius of  $P(\omega)$.
  \end{proposition}

  When $\omega=0$,  we find that $A(0)=I_n-A_0-W=I_n-V$, where $V=\sum_{i=0}^\infty A_i\sum_{j=0}^i G^j$.
  According to Theorem~4.14 in \cite{BLM}, $I-V$ is a  nonsingular M-matrix and therefore,  since $N(0)\ge 0$ and $M^{-1}\ge 0$,  $A(0)=M-N(0)$ is a regular splitting. Hence,
  the spectral radius  $\rho_0$ of $P(0)$ is less than 1. More generally, under Assumption [C1] since
  \[
  I_n -V\leq A(\omega) \leq I_n-A
  \]
  from characterization F20 in \cite{PMM} we find that  $A(\omega)$ is a nonsingular M-matrix and $A(\omega)=M-N(\omega)$ is a regular splitting. Hence,   we deduce  that  $\rho_\omega<1$. 
  The following result gives an estimate of  $\rho_\omega$ by showing its dependence as function of the relaxation parameter.

  \begin{proposition}\label{PF1}
Let $\omega$ be such that $0\le\omega\le \hat\omega$ and condition [C1] holds. Assume that the Perron eigenvector $\B v$  of 
   $P(0)$ is positive. Then we have 
   \begin{equation}\label{eq:bound}
   \rho_0-\omega (1-\rho_0)\sigma_{\max} \le \rho_\omega \le \rho_0-\omega (1-\rho_0)\sigma_{\min}
   \end{equation}
   where $\sigma_{\min}=\min_i \frac{u_i}{v_i}$ and
   $\sigma_{\max}=\max_i \frac{u_i}{v_i}$, with $\B u=(I_n-A_0)^{-1} A_1(I+G)\B v$. Moreover, $0\le \sigma_{\min},\sigma_{\max}\le \rho_0$.
   \end{proposition}
  
  \begin{proof}
   In view of the classical Collatz-Wielandt formula (see \cite{Mey}, Chapter 8) , if $P(\omega)\B v=\B w$, where $\B v>0$ and $\B w\ge 0$, then 
   \[
   \min_i \frac{w_i}{v_i}\le \rho_\omega\le \max_i \frac{w_i}{v_i}.
   \]
   Observe that
   \[
   \begin{array}{l}
   \B w=P(\omega)\B v= P(0)\B v-\omega (I_n-A_0)^{-1} A_1(I_n+G)(I-P(0))\B v=\\
   \rho_0 \B v-\omega (1-\rho_0) (I_n-A_0)^{-1} A_1(I_n+G)\B v=
   \rho_0 \B v-\omega (1-\rho_0)\B  u,
   \end{array}
   \]
   which leads to \eqref{eq:bound}, since $\B u\ge 0$. Moreover, since $A_1(I_n+G)\le W$, then
   \[
   \B u= (I_n-A_0)^{-1} A_1(I_n+G)\B v\le (I_n-A_0)^{-1} W\B v=\rho_0 \B v,
   \]
   hence $u_i/v_i\le\rho_0$ for any $i$.
  \end{proof}
  
  Observe that, in the QBD case,  where $A_i=0$ for $i\ge 2$, we have $A_1(I_n+G)= W$ and from the proof above $\B u = \rho_0\B v$.  Therefore, we have  $\sigma_{\min}=\sigma_{\max}=\rho_0$ and, hence,  $\rho_\omega=\rho_0(1-\omega(1-\rho_0))$.
  In particular, $\rho_\omega$ linearly decreases with $\omega$, and $\rho_1=\rho_0^2$.
In the general case, inequality \eqref{eq:bound} shows that the upper bound to $\rho_\omega$ linearly decreases as a function of $\omega$. Therefore the choice $\omega=\hat\omega$ gives the fastest convergence rate.

\begin{remark}\label{testr}
For the sake of illustration we consider  a  quadratic equation associated with a block tridiagonal Markov chain  taken from \cite{LR}. 
We set $A_{-1}=W+\delta I$, $A_0=A_1=W$ where $0< \delta <1$ and $W\in \mathbb R^{n\times n}$  has zero diagonal entries and  all off-diagonal entries equal  to a given  value $\alpha$ determined so that
$A_{-1}+A_0+ A_1$ is stochastic.  We find that for  $\omega_{k}=\omega \in [0, 6]$  $N(\omega)\geq 0$.   In Figure \ref{fnew} we plot  the spectral radius of $P =P(\omega)$. The linear plot  is in accordance with
Theorem \ref{PF1}. 
\begin{figure}
  \includegraphics[width=0.5\textwidth]{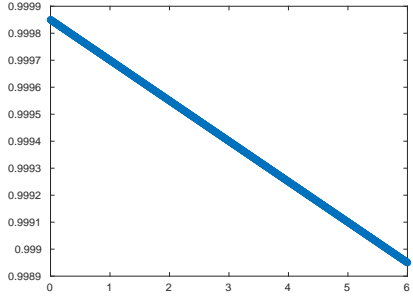}
  \caption{Plot of $\rho(P(\omega))$ for $\omega\in [0,6]$.}
\label{fnew}
\end{figure}
\end{remark}

\subsubsection{The Case $X_0$ stochastic}
In this section we analyze the convergence of the  iterative method 
\eqref{mainitrel}  starting with  a stochastic matrix $X_0$, that is, $X_0\geq 0$ and $X_0\B e=\B e$.   Eligible sequences $\{\omega_k\}_k$  are such that  $X_k\geq 0$ for any $k\geq 0$.
This  happens for $0\leq \omega_{k}\leq 1$, $k\geq 1$.  Suppose that: 
 \begin{enumerate}
  \item[S0.]\label{its0}  The sequence  $\{\omega_k\}_k$ in \eqref{mainitrel}
   is determined so that $\omega_k\geq 0$ and $X_k\geq 0$ for any $k\geq 1$.  
\end{enumerate}
Observe that the property $X_k\B e=\B e$, $k\geq 0$  is automatically satisfied. Hence,  under assumption [S0],
all the approximations generated by the iterative scheme \eqref{mainitrel} are stochastic matrices and therefore,
Proposition \ref{propsto} can be extended  in order to prove  that  the sequence $\{X_k\}_{k\in \mathbb N}$ is  convergent to $G$. 

The analysis of the  speed of convergence   follows from relation \eqref{err2}.  Let us denote 
as $\vec1(A)\in \mathbb R^{n^2}$ the vector  
obtained by stacking the columns of the matrix $A\in \mathbb R^{n\times n}$ on top of one another.  Recall that  $\vec1(ABC)=(C^T\otimes A)\vec1(B)$ for any $A,B,C \in \mathbb R^{n\times n}$.   By using this property we can rewrite \eqref{err2} as follows.  We have:
\begin{equation}\label{seq2}
\begin{array}{llllll}
[I_n \otimes (I_n-A_0)] \vec1(E_{k+1})=\omega_{k+1}[I_n\otimes A_1G(I_n-A_0)^{-1}A_1G] \vec1(E_{k})+\\
\omega_{k+1}[X_k^T\otimes A_1G(I_n-A_0)^{-1}A_1]\vec1(E_{k}) + \omega_{k+1}[Y_k^T\otimes A_1(I_n-A_0)^{-1}A_1G]\vec1(E_{k}) +\\
\omega_{k+1}[Y_k^TX_{k}^T \otimes A_1(I_n-A_0)^{-1}A_1]\vec1(E_{k})  + (1-\omega_{k+1})[I_n \otimes A_1G]\vec1(E_{k}) + \\
(1-\omega_{k+1})[X_k^T\otimes A_1]\vec1(E_{k}) + [\sum_{i=2}^\infty \sum_{j=0}^i({X_k^T}^{i-j}\otimes A_iG^j)]\vec1(E_{k}) + \\
\omega_{k+1}[\sum_{i=2}^\infty \sum_{j=0}^i({X_k^T}^{i-j}\otimes A_1G(I_n-A_0)^{-1}A_iG^j)]\vec1(E_{k}) + \\\omega_{k+1}[\sum_{i=2}^\infty \sum_{j=0}^i({Y_k^T}{X_k^T}^{i-j}\otimes A_1(I_n-A_0)^{-1}A_iG^j)]\vec1(E_{k}),
\end{array}
\end{equation}
for $k\ge 0$.
The convergence of $\{\vec1(E_k)\}$ depends on the choice of $\omega_{k+1}$, $k\geq 0$.  Suppose that $\omega_k=\omega$ for any $k\geq 0$ and [S0] holds.  Then \eqref{seq2} can be rewritten  in a compact form as
\[
\vec1(E_{k+1}) =H_k \vec1(E_{k}), \quad k\geq 0,
\]
where $H_k=H_\omega(X_k, Y_k)$ and
\[
\lim_{k\rightarrow +\infty} H_k=H_\omega(G, G)=H_\omega.
\]
It can be shown that the asymptotic rate of convergence $\sigma$ satisfies
\[
\sigma=\limsup_{k\to\infty}\left( \frac{\|\vec1(E_{k+1})\|}{ \| \vec1(E_{0}) \|}\right)^{1/k}\leq \rho(H_\omega).
\]

In the sequel we compare  the cases $\omega_k=0$, which corresponds with the  traditional fixed point iteration \eqref{eq:trad}, and
$\omega_k=1$ which reduces to the   staircase fixed point iteration \eqref{mainit}. 

For $\omega_{k+1}=0$, $k\geq 0$,  we find that 
\[
    \vec1(E_{k+1})=[\sum_{i=1}^\infty \sum_{j=0}^i({X_k^T}^{i-j}\otimes (I_n-A_0)^{-1}A_iG^j)]\vec1(E_{k}), \quad k\geq 0,  
    \]
    which means that
    \[
    H_0=\sum_{i=1}^\infty \sum_{j=0}^i({G^T}^{i-j}\otimes (I_n-A_0)^{-1}A_iG^j), \quad k\geq 0.
    \]
   Let $U^H G^TU=T$  be the Schur form of $G^T$ and  set $W= (U^H\otimes I_n)$. 
 Then 
    \[
    W \sum_{i=1}^\infty \sum_{j=0}^i({G^T}^{i-j}\otimes (I_n-A_0)^{-1}A_iG^j)
    W^{-1}=\sum_{i=1}^\infty \sum_{j=0}^i({T}^{i-j}\otimes (I_n-A_0)^{-1}A_iG^j)
    \]
    which means that  $H_0$  is similar to  the matrix
   on the right hand-side.  There follows that  the eigenvalues of $H_0$ belong  to the set 
    \[
    \cup_{\lambda}\{\mu\colon \mu \textrm{ is eigenvalue of } \sum_{i=1}^\infty \sum_{j=0}^i(\lambda^{i-j}(I_n-A_0)^{-1}A_iG^j)\}
    \]
   with $\lambda$  eigenvalue  of $G$. Since $G$ is stochastic we have 
   $|\lambda|\leq 1$. Thus,  from 
   \[
   |\sum_{i=1}^\infty \sum_{j=0}^i\lambda^{i-j}(I_n-A_0)^{-1}A_iG^j| \leq \sum_{i=1}^\infty \sum_{j=0}^i(I_n-A_0)^{-1}A_iG^j =P(0)
   \]
   we conclude that $\rho(H_0)\leq \rho(P(0))$ in view of the Wielandt theorem \cite{Mey}.

   A similar analysis can be performed in the case $\omega_k=1$, $k\geq 0$.    We find that 
\[
\begin{array}{ll} H_1=
  [I_n \otimes (I_n-A_0)^{-1}] \left\{
\sum_{i=2}^\infty \sum_{j=0}^i({G^T}^{i-j}\otimes A_iG^j) +\right.\\\left.
\sum_{i=1}^\infty \sum_{j=0}^i \left[({G^T}^{i-j}\otimes A_1G(I_n-A_0)^{-1}A_iG^j) + ({G^T}^{i-j+1}\otimes A_1(I_n-A_0)^{-1}A_iG^j)\right]
\right\}.
\end{array}
\]
By the same arguments as above we find  that the  eigenvalues of  $H_1$ belong to the set 
\[
    \cup_{\lambda}\{\mu\colon \mu \textrm{ is eigenvalue of } (I_n-A_0)^{-1} N(\lambda)\},
    \]
   with $\lambda$  eigenvalue  of $G^T$, and 
   \[
   \begin{array}{ll}N(\lambda)=\sum_{i=2}^\infty \sum_{j=0}^i(\lambda^{i-j}A_iG^j) +  \sum_{i=1}^\infty \sum_{j=0}^i(\lambda^{i-j}A_1G(I_n-A_0)^{-1}A_iG^j) + \\
    \sum_{i=1}^\infty \sum_{j=0}^i(\lambda^{i-j+1} A_1(I_n-A_0)^{-1}A_iG^j).\end{array}
   \]
   Since 
   \[
   |(I_n-A_0)^{-1} N(\lambda)|\leq P(1)
   \]
   we conclude that
   \[
   \rho(H_1)\leq \rho(P(1)).
   \]
Therefore, in the application of \eqref{mainit}  we expect a faster convergence when $X_0$ is a stochastic matrix, rather than $X_0=0$. 
Indeed, numerical results shown in Section~\ref{SEC5}   exhibit a very rapid convergence profile when $X_0$ is stochastic,  even better than the one predicted by $\rho(H_1)$. 
This might be explained with the dependence of the asymptotic convergence rate on the second eigenvalue of the
   corresponding iteration matrices as reported in  \cite{MeiB}.

\section{Adaptive  Strategies and Efficiency Analysis}\label{SEC7}

The efficiency of fixed point iterations  depends on  both speed of convergence and complexity  properties.   Markov chains  are generally   defined in terms  of  sparse matrices. To take into account this
feature we assume that $\gamma n^2$, $\gamma=\gamma(n)$,  multiplications/divisions are sufficient  to perform the following tasks:
\begin{enumerate}
\item  to compute a matrix multiplication of the form $A_i \cdot W$, where $A_i, W \in \mathbb R^{n\times n}$;
\item   to solve a linear system  of the form $(I-A_0)Z=W$, where $A_0, W \in \mathbb R^{n\times n}$.
\end{enumerate}
We also suppose that the transition matrix $P$ in \eqref{maineq-1} is banded, hence $A_i=0$ for $i>q$. This is always  the case in numerical computations where the  matrix power series   $\sum_{i=-1}^\infty A_i X_k^{i+1}$ has to be  approximated
by some  finite partial sum $\sum_{i=-1}^q A_i X_k^{i+1}$.
Under these  assumptions  we obtain the following cost estimates per step:
\begin{enumerate}
\item the traditional fixed point iteration \eqref{eq:trad} requires  $q n^3 + 2 \gamma n^2 +O(n^2)$ multiplicative operations;
  \item the $U$-based  fixed point iteration \eqref{eq:ubased} requires  $(q+4/3) n^3 +  \gamma n^2 +O(n^2)$ multiplicative operations;
 \item the staircase-based (S-based) fixed point iteration \eqref{mainit}  requires  $(q+1) n^3 +  4\gamma n^2 +O(n^2)$ multiplicative operations.
\end{enumerate}
Observe that the cost of the S-based fixed point iteration is comparable with the cost of the $U$-based iteration, which is the fastest among classical iterations \cite{BLM}. Therefore,
in the cases where the $U$/S-based  fixed point  iterative schemes require significantly less iterations to converge, these  algorithms are  more efficient than the
traditional fixed point iteration.

Concerning the relaxed versions \eqref{mainitrel} of the
S-based  fixed point iteration  for  a given fixed choice of $\omega_k=\omega$ we  get the same complexity of the  unmodified scheme
\eqref{mainit} obtained with $\omega_k=\omega=1$.  The adaptive selection of $\omega_{k+1}$ exploited in Proposition  \ref{prop2} and Remark \ref{rk2}  with  $X_0=0$ requires more care.

The strategy \eqref{strat1} is
computationally unfeasible since it needs the  additional computation of  $\sum_{i=2}^{q} A_iY_k^{i+1}$.  To approximate this quantity  we recall that
\[
A_i(Y_k^{i+1} -X_k^{i+1})=A_i\sum_{j=0}^i Y_k^j (Y_k-X_k) X_{k}^{i-j}\geq A_i\sum_{j=0}^i X_{k}^j (Y_k-X_k) X_{k-1}^{i-j}.
\]
Let $\theta_{k+1}$ be such that
\[
Y_k-X_k \geq \frac{X_k-X_{k-1}}{\theta_{k+1}}.
\] 
Then condition \eqref{strat1} can be replaced with 
\begin{equation}\label{strat1mod}
 \frac{\omega_{k+1}-1}{\omega_{k+1}} A_1((Y_k^2-X_k^2)\leq ( Y_k\Gamma_k + \Gamma_k Y_k) + (\hat \omega \theta_{k+1})^{-1} \sum_{i=2}^{q} A_i(X_{k}^{i+1}-X_{k-1}^{i+1}).
\end{equation}
The iterative scheme \eqref{mainitrel} complemented with the strategy  based on \eqref{strat1mod} for the selection of the  parameter $\omega_{k+1}$  requires no more than $(q+3) n^3 +  4\gamma n^2 +O(n^2)$
multiplicative operations.  The efficiency of this scheme will be investigated experimentally in the next section.

\section{Numerical  Results}\label{SEC8}
In this section we present the results of some numerical experiments which confirm the  effectiveness of our proposed schemes. All the algorithm have been implemented in Matlab and tested on  a PC 
i9-9900K CPU  3.60GHz$\times$8.
Our test suite includes:
\begin{enumerate}
\item Synthetic Examples:
  \begin{enumerate}
  \item The block tridiagonal Markov chain  of Remark \ref{testr}.  Observe that  the drift of the Markov chain is exactly equal to  $-\delta$.
  \item A numerical example  given in \cite{Bai} and  considered in \cite{Guo} for testing a suitable modification  --named Algorithm 1-- of the U-based fixed point iteration.
    The Markov chain of the M/G/1 type is given by
    \[
    A_{-1}=\frac{4(1-p)}{3} \left[\begin{array}{ccccc}
      0.05 & 0.1 & 0.2 & 0.3 & 0.1\\
      0.2 & 0.05 & 0.1 & 0.1 & 0.3\\
      0.1 & 0.2 & 0.3 & 0.05 & 0.1\\
      0.1 & 0.05 & 0.2 & 0.1 & 0.3\\
      0.3 & 0.1 & 0.1 & 0.2 & 0.05 \end{array}\right], \quad A_i=p A_{i-1}, ~~i\geq 0.
    \]
    The  Markov chain is positive recurrent, null recurrent or transient according as $0<p<0.5$, $p=0.5$ or $p>0.5$, respectively. In our computations we have chosen different values of $p$, in the range $0<p<0.6$ and the matrices $A_i$ are
    treated as zero matrix for $i\geq 51$.
    \item  Synthetic examples of M/G/1-type Markov chains described in
\cite{BLMfamily}. These examples are constructed in such a way that the drift of the associated Markov chain is close to a given nonnegative value.
We do not describe in detail the construction, as it would take some space, but refer the reader to~\cite[Sections~7.1]{BLMfamily}.
  \end{enumerate}
\item Application Examples:
  \begin{enumerate}
  \item   Some examples of PH/PH/1 queues  collected  in \cite[Sections~7.1]{BLMfamily}   for testing  purposes.  The  construction depends on a parameter $\rho$ with
$0\leq \rho \leq 1$. In this case the drift is $\mu=1-\rho$.
\item  The  Markov  chain  of M/G/1  type associated with the infinitesimal generator matrix $Q$ from the queuing model described in~\cite{dudin}.
  This is a complex queuing model, a BMAP/PHF/1/N model with retrial system with finite buffer of capacity $N$ and non-persistent customers.
  For the construction of the matrix $Q$ we refer the reader to~\cite[Sections~4.3 and~4.5]{dudin}.
  \end{enumerate}
  \end{enumerate}

\subsection{Synthetic Examples}

The first test concern with the validation of the analysis performed above, regarding the convergence of fixed point iterations.
In Table \ref{table1} we show  the  number of iterations required 
by  different iterative schemes  on Example  1.a with $n=100$.  Specifically  we compare the traditional fixed point iteration, the  $U$-based fixed point iteration, the S-based fixed point iteration \eqref{mainit} and
the relaxed fixed point iterations \eqref{mainitrel}.  For the latter case we consider  the  S$_\omega$-based  iteration where $\omega_{k+1}=\omega$ is fixed a priori and the
S$_{\omega(k)}$-based  iteration  where  the value of $\omega_{k+1}$  is dynamically adjusted at any step according to the strategy \eqref{strat1mod} complemented with condition  \eqref{eqom2}.
The relaxed stationary  iteration is applied for $\omega=1.8, 1.9, 2$.  The relaxed adaptive   iteration is applied with $\hat \omega=10$. The iterations are stopped when the residual error
$\| X_k-\sum_{i=-1}^q A_i X_k^{i+1}\|_\infty$ is smaller than $ tol=10^{-13}$.

\begin{table}
\begin{center}
\begin{tabular}{ |c|c|c|c|c|c|} 
\hline
$\delta$ & trad. & $U$-based & S-based & S$_\omega$-based & S$_{\omega(k)}$-based\\
\hline
\multirow{3}{*}{1.0e-2} &\multirow{3}{*}{1447}& \multirow{3}{*}{731} & \multirow{3}{*}{724} & 515 & \multirow{3}{*}{65}\\ 
&  & & &  496  &\\ 
&& & &  479  & \\
\hline
\multirow{3}{*}{1.0e-4} &\multirow{3}{*}{84067}& \multirow{3}{*}{42046} & \multirow{3}{*}{42037} & 30023 &\multirow{3}{*}{11771}\\ 
&  & & &  28976  & \\ 
&& & &  28027 & \\
\hline
\multirow{3}{*}{1.0e-6} &\multirow{3}{*}{2310370}& \multirow{3}{*}{1154998} & \multirow{3}{*}{1155352} & 825121 & \multirow{3}{*}{329843} \\ 
&  & & &  796693 & \\ 
&& & &  770250 &\\
\hline
\end{tabular}
\end{center}
\caption{Number of iterations on Example 1.a  for different values of $\delta$. The  relaxed iteration  with $\omega$ being fixed a priori is tested with
$\omega\in \{1.8, 1.9, 2\}$.}
\label{table1}
\end{table}   
The first four columns of Table 1 \ref{table1}  confirms the   theoretical comparison of asymptotic convergence rates  of classical fixed point iterations applied to a block tridiagonal matrix. Specifically,
the $U$-based and the S-based iterations are  twice faster than the traditional iteration.   Also, the relaxed stationary variants  greatly improve the convergence speed. An additional  remarkable improvement is obtained by
adjusting dynamically the value of the relaxation parameter.  Also notice that  the S$_{\omega(k)}$-based  iteration  is guaranteed to converge differently to the  stationary  S$_{\omega}$-based  iteration.

The superiority of the adaptive implementation over the other fixed point iterations is confirmed  by numerical results on Example 1.b. In Table \ref{tabnew} for different values of $p$  we show
the number of iterations  required by different iterative schemes including also Algorithm 1 implemented in  \cite{Guo}. For comparison with the results in  \cite{Guo}  here we set $tol=10^{-8}$ in the stopping criterion.

\begin{table}
\begin{center}
\begin{tabular}{ |c|c|c|c|c|c|} 
\hline
$p$ & trad. & $U$-based & S-based & Algorithm  1& S$_{\omega(k)}$-based\\
\hline
 0.3 & 14 & 11 & 10  &12& 9 \\
\hline
0.48 & 122 & 84 & 91 & 85 & 72 \\
\hline
0.5 & 7497 & 5000 & 5622 & 5001 & 4374 \\
\hline
0.55 & 53 & 37 & 39 & 39 & 32 \\
\hline
\end{tabular}
\end{center}
\caption{Number of iterations on Example 1.b  for different values of $p$.}
\label{tabnew}
\end{table}

\begin{figure}
  \centering
  \subfloat[]{\includegraphics[width=0.5\textwidth]{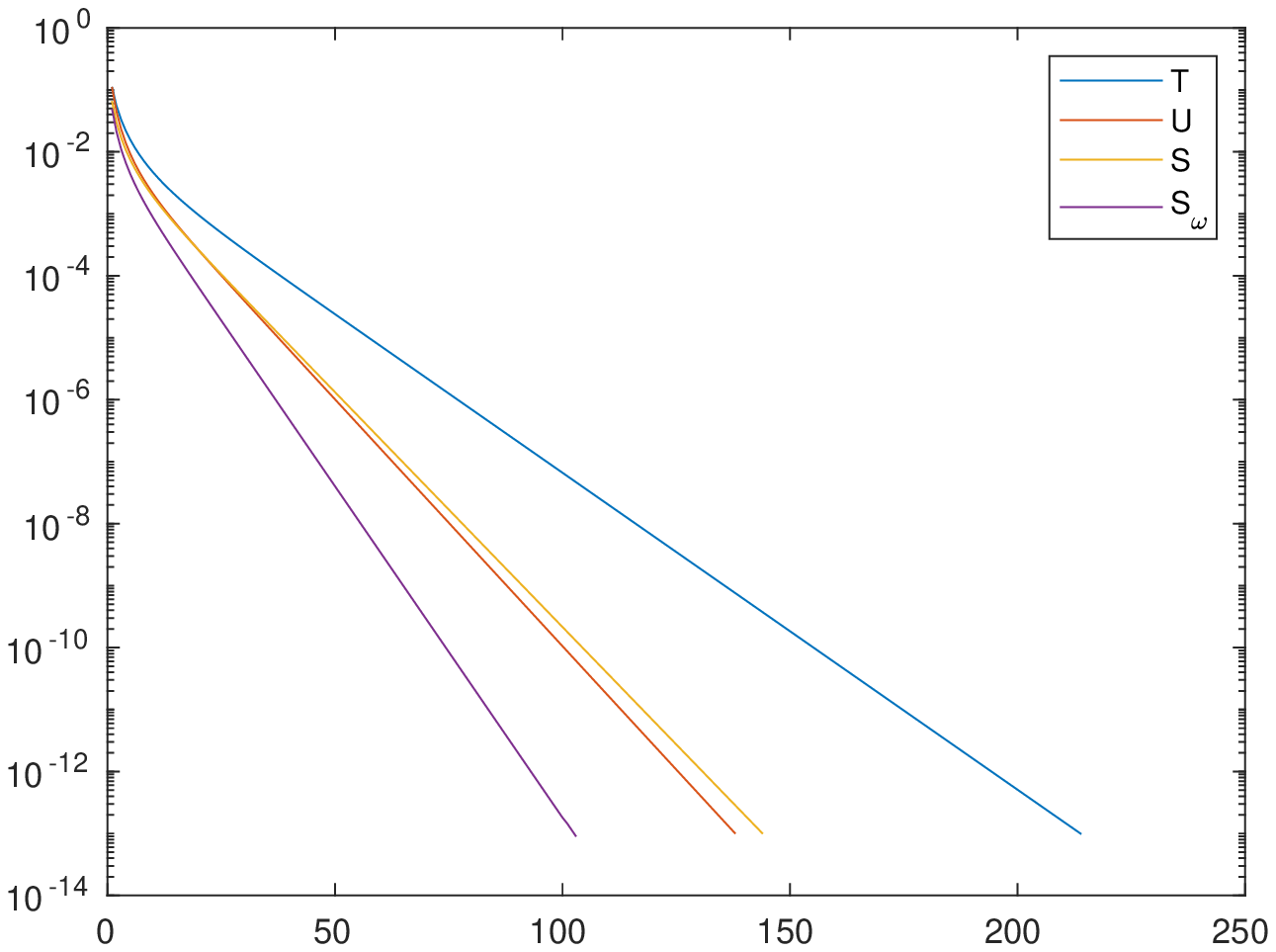}}
  \hfill
  \subfloat[]{\includegraphics[width=0.5\textwidth]{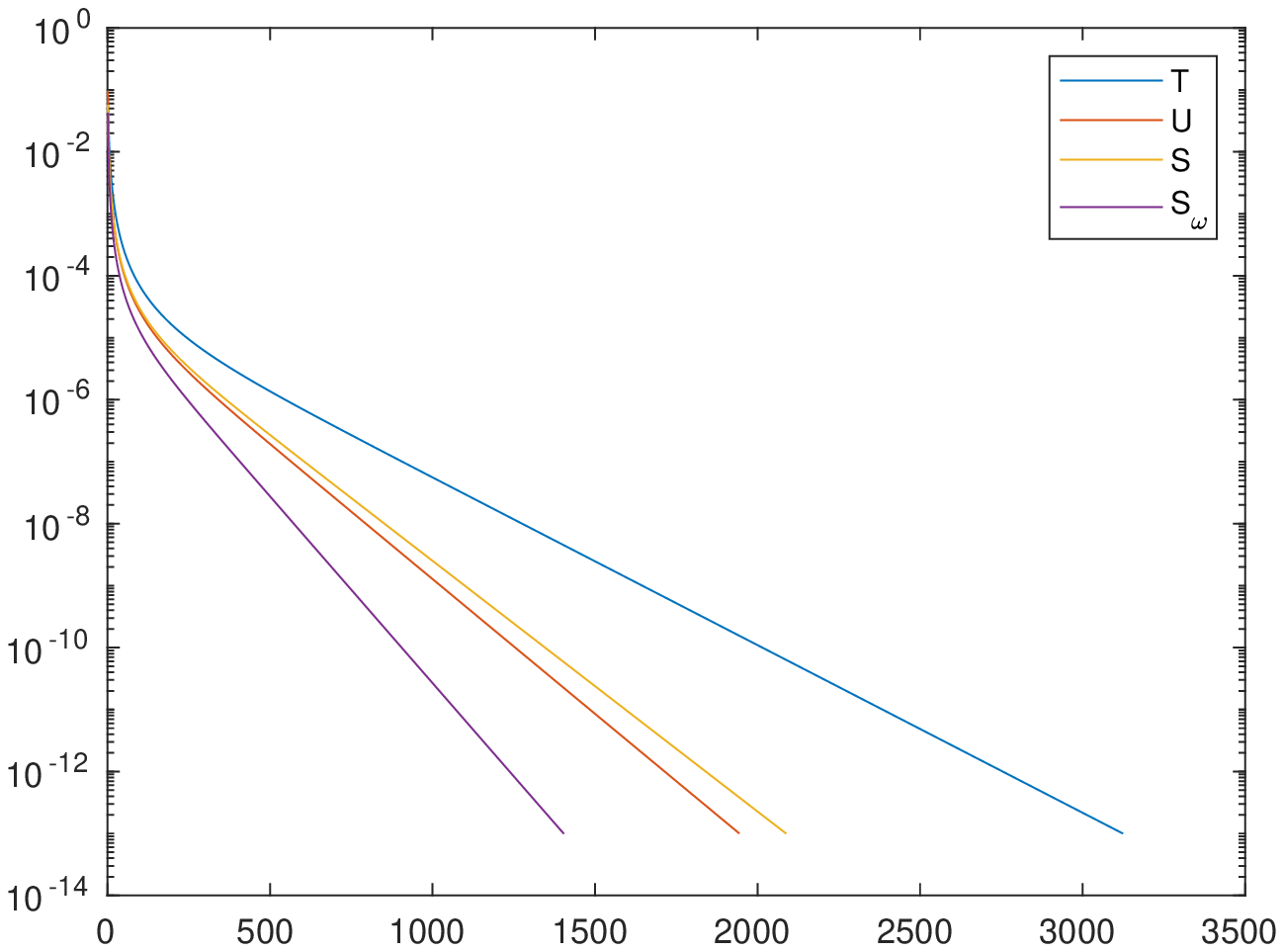}}
  \caption{Residual errors generated by the four fixed point iterations applied to the synthetic example  with drift $\mu=-0.1$ and $\mu=-0.005$.}
\label{f1}
\end{figure}

\begin{figure}
  \centering
  \subfloat[]{\includegraphics[width=0.5\textwidth]{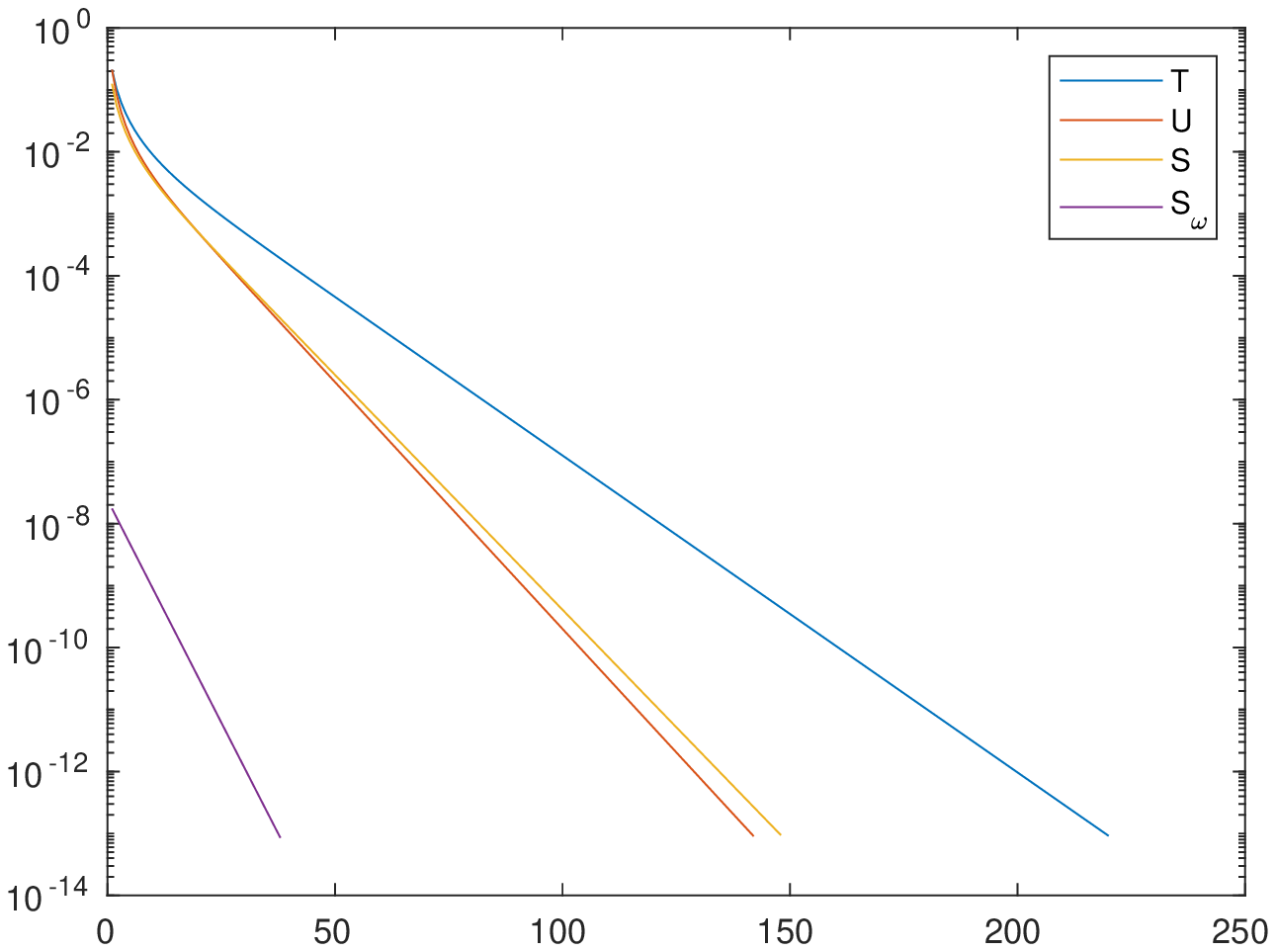}}
  \hfill
  \subfloat[]{\includegraphics[width=0.5\textwidth]{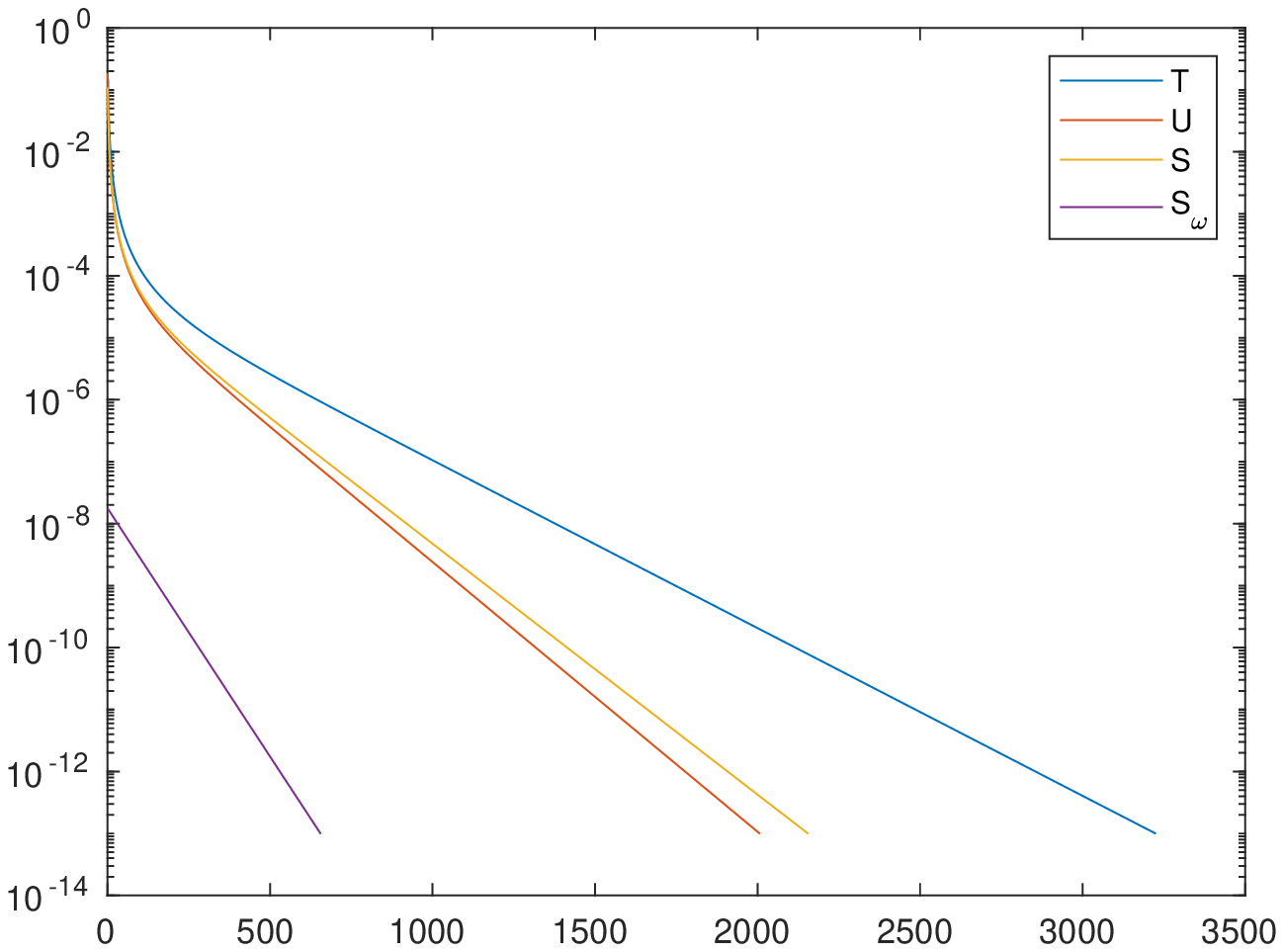}}
  \caption{Residual errors generated by the four fixed point iterations applied to the synthetic example  with drift $\mu=-0.1$ and $\mu=-0.005$.}
\label{fnew2}
\end{figure}

Finally, we  compare the convergence speed of the  traditional, $U$-based, S-based and  S$_{\omega(k)}$-based fixed point iterations  applied on  the synthetic examples of M/G/1-type Markov chains described in
\cite{BLMfamily}.  In Figure \ref{f1} we report the semilogarithmic plot of the residual error in the infinity norm  generated by the four fixed point iterations for two different values of the drift.
Observe that the adaptive relaxed  iteration is about twice faster than the traditional fixed point iteration.  The observation is confirmed in 
 in Table \ref{ts}  where we  indicate the speed-up in terms of CPU-time with respect to the traditional fixed point iteration for  different values of the drift $\mu$.
\begin{table}
\begin{center}
\begin{tabular}{ |c|c|c|c|} 
\hline
$\mu$ & $U$-based & S-based & S$_{\omega(k)}$-based\\
\hline
-0.1 & 1.6 & 1.5 & 2.1 \\
-0.05 & 1.5 & 1.4 & 2.0 \\
-0.01 & 1.6 & 1.5 & 2.1 \\
-0.005 & 1.6 & 1.5 & 2.1 \\
-0.001 & 1.6 & 1.4 & 2.2 \\
-0.0005 & 1.6 & 1.5 & 2.2 \\
-0.0001 & 1.7 & 1.6 & 2.4 \\
\hline
\end{tabular}
\end{center}
\caption{Speed-up in terms of CPU-time w.r.t. the traditional fixed point iteration for different values of the drift $\mu$. }
\label{ts}
\end{table} 
In Figure \ref{fnew2} we repeat  the set of experiments of Figure \ref{f1} with  a starting stochastic matrix $X_0= \B e\B e^T/n$.  Here the adaptive strategy is basically  the same as used before where we select  $\omega_{k+1}$  in the interval $[0, \omega_{k}]$  as  the  maximum   value which maintains the nonnegativity of $X_{k+1}$.  
In this case  the adaptive strategy seems to be quite effective in reducing the number of iterations. Other results are not so favorable and we believe that in the stochastic case the design of a general  efficient  strategy for the choice of the relaxation parameter $\omega_k$ is still an open problem. 

\subsection{Application Examples}

The first set  2.a of examples from applications includes   several cases of PH/PH/1 queues  collected   in \cite{BLMfamily}.
The  construction  of the Markov chain depends on a parameter $\rho$ with
$0\leq \rho \leq 1$  and two  integers $(i,j)$ which specify the   PH distributions of the model. The Markov chain generates in this way is denoted as Example $(i,j)$.  Its  drift is
$\mu=1-\rho$. In Tables  \ref{table11} \ref{table2} and \ref{table3} we  compare  the number of iterations for different values of $\rho$.
Here  and hereafter the relaxed stationary S$_{\omega}$-based  iteration  is applied  with  $\omega=2$. 

\begin{table}
\begin{center}
\begin{tabular}{ |c|c|c|c|c|c|} 
\hline
$\rho$ & trad. & $U$-based & S-based & S$_\omega$-based & S$_{\omega(k)}$-based\\
\hline
 0.99 & 6708 & 3666 & 3638& 1061 & 1069 \\
\hline
0.999 & 53900 & 30018 & 29271 & 8559& 8609\\
\hline
0.9999 & 386332 & 215660 & 209975& 61608& 65209\\
\hline
\end{tabular}
\end{center}
\caption{Number of iterations for different values of $\rho$ on Example $(2,6)$. }
\label{table11}
\end{table} 

\begin{table}
\begin{center}
\begin{tabular}{ |c|c|c|c|c|c|} 
\hline
$\rho$ & trad. & $U$-based & S-based & S$_\omega$-based & S$_{\omega(k)}$-based\\
\hline
 0.99 & 3735 & 3241 & 2225& 1582 & 1158 \\
\hline
0.999 & 29162 & 25241 & 17460 & 12469& 9331\\
\hline
0.9999 & 209275 & 181139 & 125689 & 89965& 68967\\
\hline
\end{tabular}
\end{center}
\caption{Number of iterations for different values of $\rho$ on Example $(8,3)$. }
\label{table2}
\end{table}

\begin{table}
\begin{center}
\begin{tabular}{ |c|c|c|c|c|c|} 
\hline
$\rho$ & trad. & $U$-based & S-based & S$_\omega$-based & S$_{\omega(k)}$-based \\
\hline
 0.99 & 11372 & 9679 & 9416& 8031& 8204\\
\hline
0.999 & 87566& 74611 & 72595 & 61992& 63346\\
\hline
0.9999 & 597330 & 509108 & 495942 & 424141& 433710\\
\hline
\end{tabular}
\end{center}
\caption{Number of iterations for different values of $\rho$ on Example $(8,9)$. }
\label{table3}
\end{table}

For comparison in Table \ref{Tsto} we show the number of iterations  on Example $(8,3)$  of Table \ref{table2}  starting with $X_0$ a stochastic matrix. We compare the traditional, $U$-based and S-based iterations.  We
observe  a rapid convergence profile and the fact that the number of iterations is independent of the drift value.

\begin{table}
\begin{center}
\begin{tabular}{ |c|c|c|c|} 
\hline
$\rho$ & trad. & $U$-based & S-based  \\
\hline
 0.99 & 52 & 45 & 31\\
\hline
0.999 & 51& 45 & 30 \\
\hline
0.9999 & 51 & 45 & 30\\
\hline
\end{tabular}
\end{center}
\caption{Number of iterations for different values of $\rho$ on Example $(8,3)$ with $X_0$ a stochastic matrix. }
\label{Tsto}
\end{table}

For a more challenging example from applications in 2.b we consider  the generator matrix $Q$ from the queuing model described in~\cite{dudin}.
In Table \ref{table4} we indicate the number of iterations for different values of the capacity $N$.

\begin{table}
\begin{center}
\begin{tabular}{ |c|c|c|c|c|c|} 
\hline
$N$ & trad. & $U$-based & S-based & S$_\omega$-based & S$_{\omega(k)}$-based \\
\hline
 20 & 4028 & 3567 & 3496& 3089 & 2826\\
\hline
30 & 4028 & 3567 & 3496& 3089 & 2826\\
\hline
40 & 4028 & 3567 & 3496& 3089& 2826 \\
\hline
\end{tabular}
\end{center}
\caption{Number of iterations for different values of $N$ on Example  described in~\cite{dudin}.  }
\label{table4}
\end{table}

\section{Conclusion and Future Work}\label{SEC9}
In this paper we have introduced a novel fixed point iteration  for solving M/G/1-type Markov chains. It is shown that this iteration complemented with suitable
adaptive relaxation techniques  is generally more efficient than other classical iterations.   Incorporating relaxation techniques into other different  inner-outer
iterative schemes as the ones introduced  in \cite{BLMfamily}  is an ongoing research. 
  
\bibliographystyle{plain} 
\bibliography{msplit}
 
\end{document}